\documentclass{article}
\usepackage{multirow}
\usepackage{amsmath,amssymb,latexsym,graphics,epsfig}
\usepackage{color}
\usepackage{amsthm}
\usepackage{graphicx,url}
\usepackage{bm}
\usepackage{hyperref}

\usepackage{tabularx}
\usepackage{caption}
\usepackage[scaled]{helvet}

\makeatletter
\DeclareTextCommand{\textcircled}{OMS}[1]{\hmode@bgroup
   \ooalign{%
      \hfil$\m@th\vcenter{\hbox{\upshape#1}}$\hfil\crcr
      \char 13 
   }%
   \vphantom{\char 13}%
\egroup}

\newtheorem{thm}{Theorem}[section]

\newtheorem{lem}[thm]{Lemma}
\newtheorem{cor}[thm]{Corollary}
\def\t{{\!\top\!}}
\def\0{{\bf 0}}
\def\1{{\bf 1}}
\def\x{{\bf x}}
\def\y{{\bf y}}
\def\v{{\bf v}}
\def\e{{\bf e}}
\def\F2{\mathbb{F}_{2}}
\def\rank{\mbox{\rm 2-rank}}
\def\col{{\mbox{\rm Col$_2$}}}
\def\H{\otimes}       
\def\qed{\hfill$\Box$\\}

\begin{document}

\title
{\bf Graph switching, 2-ranks, and \\ graphical Hadamard matrices}
\author{Aida Abiad$^a$, Steve Butler$^b$ and Willem H. Haemers$^c$
\\ \\
{\small $^a$Department of Quantitative Economics} \\
{\small Maastricht University, Maastricht, The Netherlands}\\
{\small {\tt
aidaabiad@gmail.com}} \\
{\small $^b$Department of Mathematics} \\
{\small Iowa State University, Ames, IA, USA}\\
{\small {\tt
butler@iastate.edu}} \\
{\small $^c$Department of Econometrics and Operations Research} \\
{\small Tilburg University, Tilburg, The Netherlands}\\
{\small {\tt
haemers@uvt.nl}}\\
}

\date{}
\maketitle

\begin{abstract}
\noindent
We study the behaviour of the 2-rank of the adjacency matrix of a graph
under Seidel and Godsil-McKay switching, and apply the result to graphs coming from graphical Hadamard matrices of order $4^m$.
Starting with graphs from known Hadamard matrices of order $64$, we find (by computer) many Godsil-McKay switching sets
that increase the 2-rank.
Thus we find strongly regular graphs with parameters $(63,32,16,16)$, $(64,36,20,20)$, and $(64,28,12,12)$ for almost all feasible 2-ranks.
In addition we work out the behaviour of the 2-rank for a graph product related to
the Kronecker product for Hadamard matrices, which enables us to find many graphical
Hadamard matrices of order $4^m$ for which the related strongly regular graphs
have an unbounded number of different 2-ranks.
The paper extends results from the article {\lq}Switched symplectic graphs and their 2-ranks{\rq} by the first and the last author.
\\[5pt] \noindent
\end{abstract}

\section{Introduction}

The \emph{2-rank} of a graph is the rank of its adjacency matrix over $\F2$.
It is a well-studied and useful graph parameter (see for example \cite{BE,P}).
Sometimes the 2-rank can be used to distinguish cospectral graphs, such as
strongly regular graphs (for short SRGs) with the same parameters (and therefore the same spectrum).
An important fact is that the 2-rank of a graph is an even number (see \cite{BE}, or \cite{HPR}).

Godsil-McKay switching (for short GM-switching) is an operation on graphs that does not change the spectrum of the adjacency matrix.
For GM-switching to work, one needs a vertex subset with special properties, called a GM-set.
However, GM-switching can change the 2-rank, in which case the switched graph is obviously non-isomorphic to the original one.
This idea was a starting point of an earlier paper~\cite{AH} by two of the present authors.
They gave switching sets in the symplectic graph $Sp(2m,2)$, which is a famous SRG with parameters
\[
P_0(m)=(2^{2m}-1,\ 2^{2m-1},\ 2^{2m-2},\ 2^{2m-2}),
\]
which increase the 2-rank after switching.
In addition, repeated GM-switching was applied for the case $m=3$, and many new strongly regular graphs
with parameters $P_0(3)=(63,32,16,16)$ were found and the 2-ranks vary from 6 to 18.
In this paper we use an improved computer search and obtain examples with 2-rank 20, 22 and 24.
In addition we apply the same idea to SRGs with parameters
\[
P_{\pm}(m)= (2^{2m},\ 2^{2m-1}\pm 2^{m-1},\ 2^{2m-2}\pm 2^{m-1},\ 2^{2m-2}\pm 2^{m-1}).
\]
For $m=3$ we found such SRGs for all 2-ranks in $\{8,10,\ldots,26\}$.
SRGs with parameter sets $P_0(m)$ and $P_\pm(m)$ correspond to graphical Hadamard matrices of order $4^m$.
For these Hadamard matrices there is a recursive construction using Kronecker products.
We find the behaviour of the 2-rank of the corresponding graphs for this construction.
Using this we obtain SRGs with parameters $P_0(m)$ and 2-ranks $2m,2m+2,\ldots,2m+18\lfloor m/3\rfloor $ and
SRGs with parameters $P_\pm(m)$ and 2-ranks $2m+2,2m+4,\ldots,2m+2+18\lfloor m/3\rfloor$.
It is known that the $2$-ranks of SRGs with parameter sets $P_0(m)$ and $P_\pm(m)$ lie in the intervals
$\left[2m\ ,\ 2^{2m-1}-2^{m-1}-2\right]$, and $\left[2m+2\ ,\ 2^{2m-1}-2^{m-1}\right]$, respectively (see \cite{HPR} and \cite{AH}).
For $m=2$, the upper and lower bound coincide, and for $m=3$ there are ten possible $2$-ranks for each parameter set,
of which only one value is still open ($26$ for $P_0(3)$, and $28$ for $P_\pm(3)$).

For the relevant background on graphs and matrices we refer to~\cite{BH}.
The $m\times n$ all-ones matrix is denoted by $J_{m,n}$, or just $J$, and $\1$ is the all-ones vector.
We denote the column space of a matrix $M$ over $\F2$ by $\col(M)$.
If $G$ is a graph with adjacency matrix $A$, then we sometimes write $\col(G)$ instead of $\col(A)$.

\section{Seidel switching}

Consider a graph $G=(V,E)$ of order $n$ and let $X$ be a subset of $V$ of cardinality $m$ ($0<m<n$).
{\em Seidel switching} in $G$ with respect to $X$ is an operation on $E$ defined as follows:
All edges from $E$ between $X$ and $V\setminus X$ are deleted, and all possible edges between $X$ and $V\setminus X$
which are not in $E$ are inserted (edges with both vertices inside $X$, or outside $X$ remain unchanged).
If $A$ is the adjacency matrix of $G$, then $S=J-2A-I$ is the {\em Seidel matrix} of $G$.
So the off-diagonal entries of $S$ are $\pm 1$, and $S_{i,j}=-1$ if and only if $i$ and $j$ are adjacent.
In terms of the Seidel matrix,
Seidel switching with respect to $X$ means that the rows and columns corresponding to $X$ are multiplied by $-1$.
This implies that Seidel switching does not change the spectrum of the Seidel matrix $S$.

Assume that the subset $X$ corresponds to the first rows and columns of $A$,
and let $A_X$ denote the adjacency matrix of the switched graph $G_X$.
Then
\[
A_X = A+K\ (\mbox{mod}\ 2),\mbox{ where } K=\left[\begin{array}{cc} O & J_{m,n-m}\\J_{n-m,m} & O \end{array}\right].
\]
We know that $\rank(A)$ is even, and since rank$(K)=2$ (over any field),
it follows that $\rank(A_X)\in\{\rank(A)-2,\ \rank(A),\ \rank(A)+2\}$.

\begin{lem}\label{Seidel2}
Suppose $G_X$ is obtained from $G$ by Seidel switching with respect to the set $X$ of neighbors of a vertex $x$ of $G$.
Then $x$ is an isolated vertex of $G_X$, and $\rank(G_X) = \rank(G)-2$ if $\1\in\col(G)$, and $\rank(G_X) = \rank(G)$ otherwise.
\end{lem}

\noindent
{\bf Proof.}
The first claim is obvious.
Let $\x$ be the characteristic vector of $X$.
Then clearly $\x\in\col(A)$, and since $\1$ and $\x$ span $\col(K)$, we have
\[
\col([A_X\ \1\ \x]) = \col([A\ \1\ \x]) = \col([A\ \1]).
\]
Suppose $\1\in\col(A)$.
Then $\col(A) = \col([A\ \1]) = \col([A_X\ \1\ \x])$.
The switched graph $G_X$ has an isolated vertex, therefore $\1\notin\col(A_X)$.
Hence $\col(A_X)$ is a proper subspace of $\col(A)$,
from which it follows that $\rank(A)=\rank(A_X)+2$.

If $\rank(A_X) = \rank(A)-2$, then $\col([A_X\ \1\ \x]) = \col([A\ \1\ \x])$
implies that $\x, \1 \not\in \col(A_X)$ and $\x,\1\in\col(A)$.
\qed

\section{Godsil-McKay switching}\label{GM}
Godsil and McKay introduced the following switching operation that leaves the spectrum of the adjacency matrix invariant.
\begin{thm}\label{GMswitching}
Let $G$ be a graph and let $W$ be a subset of the vertex set of $G$ which induces a regular subgraph.
Assume that each vertex outside $W$ is adjacent to $|W|$, $\frac{1}{2}|W|$ or $0$ vertices of $W$.
Make a new graph $G_W$ from $G$ as follows.
For each vertex $v$ outside $W$ with $\frac{1}{2}|W|$ neighbors in $W$,
delete the $\frac{1}{2}|W|$ edges between $v$ and $W$, and join $v$ instead to the $\frac{1}{2}|W|$
other vertices in $W$.
Then $G$ and $G_W$ have the same adjacency spectrum.
\end{thm}
The operation that changes $G$ into $G_W$ is called \emph{Godsil-McKay switching} (for short GM-switching).
Notice that if all vertices outside $W$ have $\frac{1}{2}|W|$ neighbors in $W$, then GM-switching is a special case of Seidel switching.

It is well-known that if a graph $G_W$ has the same spectrum as a SRG $G$,
then $G_W$ is also strongly regular with the same parameters as $G$. 
Therefore GM-switching provides a tool to construct new SRGs from known ones.
However, $G_W$ may be isomorphic with $G$,
but if GM-switching changes the 2-rank, this is obviously not the case.

Similar to Seidel switching, GM-switching can be described in terms of the adjacency matrices $A$ and $A_W$
of $G$ and $G_W$.
Indeed, $A_W=A+L$ (mod~2), where $L$ is defined by
$L_{i,j}=1$ if $i\in W$, $j\not\in W$ and $j$ has $\frac{1}{2}|W|$ neighbors in $W$; otherwise $L_{i,j}=0$.
Then $\rank(L)=2$ and therefore
$\rank(G_W) \in \{\rank(G)-2,\ \rank(G),\ \rank(G)+2\}$ (see also~\cite{AH}).
Moreover, in the case $\rank(G)$ increases, we have $\col(A)\subset\col(A_W)$,
and therefore $\1\in\col(A)$ implies $\1\in\col(A_W)$.

\section{Hadamard matrices}\label{Hadamard}
A square $(+1,-1)$-matrix $H$ of order $n$ is a \emph{Hadamard matrix} whenever $HH^{\t}=nI$.
If a row or a column of a Hadamard matrix is multiplied by $-1$, it remains a Hadamard matrix.
We can apply this operation a number of times such that the first row and column consist of all ones.
Such a Hadamard matrix is called \emph{normalized}.
A Hadamard matrix $H$ is said to be \emph{graphical} if $H$ is symmetric and it has constant diagonal,
and $H$ is {\em regular} if all row and column sums are equal.
We assume that the diagonal entries of a graphical Hadamard matrix $H$ are equal to $1$ (otherwise consider $-H$).
Then $A_H=\frac{1}{2}(J-H)$ is the adjacency matrix of a graph, say $G_H$.
Note that $H-I$ is the Seidel matrix of $G_H$.
If $H$ is normalized, then $G_H$ has an isolated vertex, and it is well-known that for $n>4$
the graph on the remaining $n-1$ vertices is strongly regular with parameters
$(n-1,\frac{n}{2},\frac{n}{4},\frac{n}{4})$.
If $H$ is graphical and regular, then the row and column sums are equal to $\epsilon\sqrt{n}$ where $\epsilon=\pm 1$,
and $G_H$ is strongly regular graph with parameters
$(n,\frac{n}{2}-\frac{\epsilon}{2}\sqrt{n},\frac{n}{4}-\frac{\epsilon}{2}\sqrt{n},\frac{n}{4}-\frac{\epsilon}{2}\sqrt{n})$.
Conversely, any strongly regular graph with one of the above parameters comes from a Hadamard matrix in the described way.

It is well known that if $H_1$ and $H_2$ are Hadamard matrices, then so is the Kronecker product $H_1\otimes H_2$.
Moreover, if $H_1$ and $H_2$ are normalized, then so is $H_1\otimes H_2$, if $H_1$ and $H_2$ are graphical,
then so is $H_1\otimes H_2$, and if $H_1$ and $H_2$ are regular then so is $H_1\otimes H_2$.
For example
\[
H_1=
{\scriptsize
\left[\begin{array}{rrrr}
1 &\! -1 & 1 & 1 \\
-1 & 1 & 1 & 1 \\
1 & 1 & 1 & \!-1 \\
1 & 1 & \!-1 & 1
\end{array} \right]
}
\mbox{ and }
H_2=
{\scriptsize
\left[\begin{array}{rrrr}
 1 &\! -1 &\! -1 &\! -1 \\
\! -1 & 1 &\! -1 &\! -1 \\
\! -1 &\! -1 & 1 &\! -1 \\
\! -1 &\! -1 &\! -1 & 1
\end{array} \right]
}
\]
are regular graphical Hadamard matrices, and so are $H_1\otimes H_1$, $H_1\otimes H_2$, and $H_2\otimes H_2$.
The SRGs $G_{H_1\otimes H_1}$ and $G_{H_2\otimes H_2}$ are isomorphic with parameters $P_-(2)$.
The graph is known as the lattice graph $L(4)$.
The SRG $G_{H_1\otimes H_2}$ has parameters $P_+(2)$, and is known as the Clebsch graph.
For later use we define $G_-(3)=G_{H_1\otimes H_1\otimes H_1}$, and $G_+(3)=G_{H_1\otimes H_1\otimes H_2}$,
which are SRGs with parameters $P_-(3)$ and $P_+(3)$, respectively.

For a recent survey on graphical Hadamard matrices, we refer to~\cite{B}.

\section{A graph product and its 2-rank behaviour}\label{prod}
Inspired by the Kronecker product for Hadamard matrices we define the graph product denoted by $\H$ as follows.
For $i=1,2$ let $G_i$ be a graph of order $n_i$ with vertex set $V_i$, Seidel matrix $S_i$ and adjacency matrix $A_i$.
Then $G_1\H G_2$ is the graph with vertex set $V_1\times V_2$, where two vertices $(x_1,x_2)$ and $(y_1,y_2)$ are adjacent
whenever $\{x_i,y_i\}$ is an edge in $G_i$ for $i=1,2$, or  $\{x_i,y_i\}$ is not an edge in $G_i$ for $i=1,2$.
Thus the Seidel matrix of $G_1\H G_2$ equals $(S_1+I)\otimes(S_2+I) - I$.
So if $H_1$ and $H_2$ are graphical Hadamard matrices, then $G_{H_1}\H G_{H_2}=G_{H_1\otimes H_2}$.

\begin{thm}\label{prod2rank}
For two graphs $G_1$ and $G_2$ the following hold:
\begin{description}
\item[(i)] $\1\in \col(G_1\H G_2)$ if and only if $\1\in \col(G_1)$ or $\1\in\col(G_2)$,
\item[(ii)] if $\1 \in\col(G_1)$ and $\1\in\col(G_2)$ then
\[
\rank(G_1\H G_2)=\rank(G_1)+\rank(G_2)-2,
\]
\item[(iii)] if $\1\not\in\col(G_1)$ or $\1\not\in\col(G_2)$ then
\[
\rank(G_1\H G_2)=\rank(G_1)+\rank(G_2).
\]
\end{description}
\end{thm}
%

\noindent
{\bf Proof.}
Let $n_i$ be the number of vertices of $G_i$ for $i=1,2$,
and let $A_1$, $A_2$, and $A_{1,2}$ be the adjacency matrix of $G_1$, $G_2$ and $G_1\H G_2$, respectively.
Then over $\F2$ the matrix $A_{1,2}$ satisfies
\begin{eqnarray}\label{Hprod}
A_{1,2} = A_1\otimes J_{n_2,n_2} + J_{n_1,n_1}\otimes A_2.
\end{eqnarray}
(i)~Assume $\1\in\col(A_1)$, then $A_1\v=\1$ for some $\v$ in $\F2^{n_1}$.
The weight of $\v$ is equal to $\1^\top\v=\v^\top A_1 \v=0$~(mod~$2$),
because $A_1$ is symmetric with zero diagonal.
If $\e$ is a unit vector, and $\v'=\v\otimes\e$, then (\ref{Hprod}) implies that (over $\F2$)
\[
A_{1,2}\v'=(A_1\otimes J)(\v\otimes\e) + (J\otimes A_2)(\v\otimes\e) = A_1\v\otimes \1 + J\v\otimes A_2\e = \1+\0.
\]
Therefore $\1\in\col(A_{1,2})$.
Conversely, assume $\1\in\col(A_{1,2})$.
Then (\ref{Hprod}) implies that there exist $\v_1\in\F2^{n_1}$ and $\v_2\in\F2^{n_2}$
such that $\1=A_1\v_1\otimes\1 + \1\otimes A_2\v_2$.
Therefore $A_i\v_i=\alpha_i\1$ with $\alpha_i\in\F2$ for $i=1,2$.
Clearly $\alpha_1$ or $\alpha_2$ is nonzero, so $\1\in\col(A_1)$ or $\1\in\col(A_2)$.

To prove (ii) and (iii), we first assume that $G_1$ and $G_2$ both have an isolated vertex.
Then clearly $\1\not\in\col(A_1)$ and $\1\notin\col(A_2)$.
For $i=1,2$, let $r_i$ be $\rank(A_i)$, and let $V_i$ be a $n_i\times r_i$ submatrix of $A_i$,
such that its columns are a basis for $\col(A_i)$.
Consider the matrix
\[
V_{1,2} = \left[
\begin{array}{c|c}
V_1\otimes J_{n_2,r_2} & J_{n_1,r_1}\otimes V_2
\end{array}\right].
\]
Since $A_1$ and $A_2$ have a zero column, the columns of $V_{1,2}$ are columns of $A_{1,2}$,
and by (\ref{Hprod}) they span $\col(A_{1,2})$.
Also the columns of $V_{1,2}$ are independent, since $ \col(V_1\otimes J_{n_2,r_2})$ and $\col(J_{n_1,r_1}\otimes V_2)$
have no nonzero vector in common.
Therefore $\rank(A_{1,2})=r_1+r_2$.

If $G_1$ or $G_2$ has no isolated vertex, we apply Seidel switching.
Suppose that for $i=1,2$ $G_i'$ is obtained from $G_i$ by Seidel switching with respect to the neighbors of a vertex $x_i$.
Then $x_i$ is an isolated vertex of $G_i'$, and it follows straightforwardly that
$G'_1\H G'_2=(G_1\H G_2)'$, where $(G_1\H G_2)'$ is obtained from $G_1\H G_2$ by Seidel switching with respect to the neighbors of
$(x_1,x_2)$.
Now we use Lemma~\ref{Seidel2}.
If $\1\in\col(A_1)$ and $\1\in\col(A_2)$ then $\1\in \col(A_{1,2})$, $\rank(G_1')=\rank(G_1)-2$, $\rank(G_2')=\rank(G_2)-2$
and $\rank(G_1\H G_2)=\rank(G_1\H G_2)'+2$.
Therefore $\rank(G_1\H G_2) = \rank(G_1)+\rank(G_2)-2$, which proves (ii).
The cases of statement (iii) go similarly.
\qed

\section{SRGs with parameters $P_0(3)$ and $P_\pm(3)$}

In this section, we report the result of a computer search for GM-switching sets in SRGs with parameters
$P_0(3)=(63,32,16,16)$, $P_+(3)=(64,36,20,20)$, and $P_-(3)=(64,28,12,12)$.
We start with known SRGs with the smallest possible 2-rank and search for GM-switching sets of size 4
that increase the 2-rank after switching.
We switch, and then continue the search with the newly obtained SRGs.
However, unlike in the preceding paper~\cite{AH}, we do not stop if we find no switching set that increases the 2-rank.
Instead, we also consider switching sets that do not change the 2-rank, switch and then continue the search.
A complete search considering all suitable switching sets of size 4 in each step is far out of reach,
so we stop the search if we have not found a switching set that increases the 2-rank in several thousand iterations.

For more details about the computational aspects, see the SAGE
worksheet\footnote{\url{https://cocalc.com/projects/57b6e497-d392-406c-aa9c-80221136762e/files}},
where graph strings and series of switching sets (following SAGE vertex labelling) are provided in order to reproduce
the results shown in this section. 

\begin{table}
{\small
\[
\begin{array}{cc}
\mbox{GM-switching set} & \rank\\ 
\hline
\\[-5pt]
\{(100000),\ (010000),\ (101000),\ (011000)\}&  {\ 8}\\
\{(000010),\ (000001),\ (001010),\ (001001)\}&  {10} \\
\{(100010),\ (101010),\ (110011),\ (111011)\}&  {12}\\
\{(000100),\ (010100),\ (001111),\ (011111)\}&  {14}\\
\{(000110),\ (000101),\ (010110),\ (010101)\}&  {16}\\
\{(001000),\ (100001),\ (110010),\ (011011)\}&  {18}\\
\{(110100),\ (111100),\ (100111),\ (101111)\}&  {18}\\
\{(110100),\ (111100),\ (110101),\ (111101)\}&  {20}\\
\{(010100),\ (110110),\ (101101),\ (001111)\}&  {20}\\
\{(100100),\ (110100),\ (101100),\ (111100)\}&  {22} \\
\{(000011),\ (110001),\ (001011),\ (111001)\}&  {22}\\
\{(000001),\ (001001),\ (110001),\ (111001)\}&  {22}\\
\{(010000),\ (000001),\ (010010),\ (000011)\}&  {24}
\end{array}
\]
}
\caption{Increasing 2-ranks by repeated GM-switching in $Sp(6,2)$}\label{P0}
\end{table}

A SRG with parameters $P_0(3)$ has a minimal possible 2-rank of 6 and there is a unique such SRG (see \cite{P}):
the symplectic graph $Sp(6,2)$.
The vertex set $V$ of $Sp(6,2)$ consist of the nonzero vectors in $\F2^{6}$,
and two vertices $\x=(x_1,\ldots,x_{6})$ and $\y=(y_1,\ldots,y_{6})$
are adjacent if $x_1y_2+x_2y_1+x_3y_4+x_4y_3+x_{5}y_{6}+x_{6}y_{5}=1$.
In Table~\ref{P0}, the first row gives a GM-switching set in $Sp(6,2)$,
and each subsequent row gives a GM-switching set in the SRG corresponding to the resulting graph from carrying out GM-switching on the previous row.
The last column gives the $\rank$ after switching.
Note that at some stages we use switching sets that do not increase the $\rank$.
Here the upper bound for the 2-rank is 26.
Unfortunately our search found no such graph, 
so the existence of a SRG with $\rank$ $26$ and parameters $P_0(3)$ remains open.

We know two nonisomorphic SRGs with parameters $P_-(3)$ and $\rank$ 8.
One is $G_-(3)=2K_2 \H 2K_2 \H 2K_2$, which was defined in Section~\ref{Hadamard}.
We easily have $\rank(2K_2)=4$, and $\1\in\col(2K_2)$, so Theorem~\ref{prod2rank}(ii) gives $\rank(G_-(3))=8$.
Let $\{1,2,3,4\}$ be the vertex set of $2K_2$, and let $\{1,2\}$ and $\{3,4\}$ be the edges.
Then each vertex of $G_-(3)$ can be represented by a triple in $\{1,2,3,4\}^3$.
With this notation the GM-switching sets that lead to a SRG with parameter set $P_-(3)$ and $\rank$ 26
are given in the left part of Table~\ref{P-}.

\begin{table}
{\small
\[
\hspace{-8pt}
\begin{array}{l|l}
\begin{array}{cc}
\mbox{GM-switching sets} & \!\!\mbox{2-rk}\!\!\!\\
\hline\\[-5pt]
\{(1, 1, 1), (1, 1, 3), (2, 2, 1), (2, 2, 3) \}&{10}\\
\{(1, 1, 1), (1, 1, 3), (2, 1, 1), (2, 1, 3) \}&{12}\\
\{(1, 1, 2), (2, 2, 4), (4, 2, 1), (3, 1, 3) \}&{14}\\
\{(4, 4, 4), (1, 2, 2), (2, 3, 2), (3, 1, 4) \}&{16}\\
\{(4, 3, 4), (3, 3, 2), (3, 2, 4), (4, 2, 2) \}&{18}\\
\{(4, 4, 2), (3, 4, 2), (3, 1, 1), (4, 1, 1) \}&{20}\\
\{(1, 1, 2), (1, 1, 4), (2, 2, 2), (2, 2, 4) \}&{20}\\
\{(4, 4, 1), (3, 3, 3), (1, 3, 2), (2, 4, 4) \}&{22}\\
\{(1, 3, 2), (1, 3, 4), (2, 4, 2), (2, 4, 4) \}&{22}\\
\{(4, 3, 3), (3, 3, 1), (4, 2, 3), (3, 2, 1) \}&{24}\\
\{(2, 1, 2), (4, 3, 2), (3, 2, 2), (1, 4, 2) \}&{24}\\
\{(3, 4, 4), (1, 2, 4), (2, 3, 4), (4, 1, 4) \}&{26}
\end{array}
&
\begin{array}{cc}
\mbox{GM-switching sets} & \!\!\mbox{2-rk}\!\!\!\\
\hline\\[-5pt]
\{(1,1,1), (3,2,2), (4,4,1), (2,3,2)\}&{10}\\
\{(3,3,1), (2,2,1), (3,3,2), (2,2,2)\}&{12}\\
\{(2,2,4), (3,3,4), (2,3,4), (3,2,4)\}&{14}\\
\{(3,2,2), (2,3,2), (4,4,4), (1,1,4)\}&{16}\\
\{(1,4,2), (1,1,2), (4,1,3), (4,4,3)\}&{18}\\
\{(4,2,2), (2,4,2), (3,1,3), (1,3,3)\}&{20}\\
\{(1,3,2), (4,3,2), (4,2,3), (1,2,3)\}&{20}\\
\{(1,2,1), (2,4,1), (4,3,4), (3,1,4)\}&{22}\\
\{(4,3,1), (2,4,1), (1,2,4), (3,1,4)\}&{22}\\
\{(4,2,1), (3,4,2), (1,3,1), (2,1,2)\}&{24}\\
\{(1,4,1), (4,1,1), (2,2,4), (3,3,4)\}&{24}\\
\{(3,2,2), (3,3,1), (2,2,4), (2,3,3)\}&{26}
\end{array}
\end{array}
\]
}
\caption{Increasing 2-ranks by repeated GM-switching in $G_-(3)$ (left) and $G'_-(3)$ (right)}\label{P-}
\end{table}
As mentioned before, $2K_2\H 2K_2$ is an SRG with parameters $P_-(2)$ known as the lattice graph $L(4)$.
However there is one other SRG with parameters $P_-(2)$,
known as the Shrikhande graph, which can be obtained from $L(4)$ by Seidel switching with respect to any
$4$-coclique (in this particular case, Seidel switching and GM-switching are the same).
We easily have $\rank({\rm Shrikhande})=6$ and $\1\in\col({\rm Shrikhande})$.
Define $G'_-(3)=\mbox{Shrikhande}\H 2K_2$.
Then $G'_-(3)$ is another SRG with parameters $P_-(3)$, $\rank(G'_-(3))=8$ and $\1\in\col(G'_-(3))$.
We have also searched for GM-switching sets in $G'_-(3)$.
The outcome is given in the right part of Table~\ref{P-},
where we use the same vertex set as for $G_-(3)$, but replaced $2K_2\H 2K_2$
by the Shrikhande graph obtained by switching with respect to $\{(1,1),(2,2),(3,3),(4,4)\}$.

We also considered two nonisomorphic SRGs with parameter sets $P_+(3)$ and 2-rank 8.
The first one is $G_+(3)=2K_2\H 2K_2\H K_4$, which was defined in Section~\ref{Hadamard}.
The other one is $G'_+(3)= \mbox{Shrikhande}\H K_4$.
Again the vertex set is given by $\{1,2,3,4\}^3$.
The sequence of GM-switching sets leading to SRGs with parameters $P_+(3)$ and $\rank$ 26 is given in Table~\ref{P+}.

\begin{table}
{\small
\[
\hspace{-8pt}
\begin{array}{l|l}
\begin{array}{cc}
\mbox{GM-switching sets} & \!\!\mbox{2-rk}\!\!\!\\
\hline\\[-5pt]
\{(1,1,1),(1,1,3),(1,2,1),(1,2,3)\}&{10}\\
\{(2,1,3),(1,3,1),(3,4,3),(4,2,1)\}&{12}\\
\{(2,1,1),(3,3,4),(2,4,2),(3,2,3)\}&{14}\\
\{(2,1,2),(2,1,4),(4,4,2),(4,4,4)\}&{16}\\
\{(2,3,3),(4,3,2),(1,4,4),(3,4,1)\}&{18}\\
\{(3,4,4),(4,4,3),(2,4,1),(1,4,2)\}&{20}\\
\{(3,1,2),(4,3,2),(1,4,4),(2,2,4)\}&{20}\\
\{(1,1,4),(4,3,1),(1,4,3),(4,2,2)\}&{20}\\
\{(4,1,1),(3,3,3),(1,4,1),(2,2,3)\}&{22}\\
\{(4,1,3),(3,1,1),(4,3,3),(3,3,1)\}&{22}\\
\{(2,1,1),(4,1,4),(3,2,3),(1,2,2)\}&{24}\\
\{(2,3,3),(4,3,2),(1,4,4),(3,4,1)\}&{24}\\
\{(1,3,3),(3,3,2),(2,4,4),(4,4,1)\}&{24}\\
\{(2,4,4),(3,4,3),(4,2,1),(3,2,3)\}&{24}\\
\{(2,1,2),(2,1,4),(4,3,4),(1,4,4)\}&{26}
\end{array}
&
\begin{array}{cc}
\mbox{GM-switching sets} & \!\!\mbox{2-rk}\!\!\!\\
\hline\\[-5pt]
\{(1,1,1),(2,1,2),(4,4,4),(3,4,3)\}&{10}\\
\{(1,1,2),(2,1,1),(2,2,2),(1,2,1)\}&{12}\\
\{(1,2,3),(4,4,1),(1,4,1),(4,2,3)\}&{14}\\
\{(1,2,4),(3,3,2),(1,1,3),(3,4,1)\}&{16}\\
\{(3,3,1),(2,1,3),(3,2,2),(2,4,4)\}&{18}\\
\{(2,2,3),(4,4,2),(2,3,4),(4,1,1)\}&{20}\\
\{(3,1,3),(4,2,1),(1,4,3),(2,3,1)\}&{20}\\
\{(2,4,1),(4,2,2),(1,4,4),(3,2,3)\}&{22}\\
\{(2,2,4),(3,4,2),(4,2,1),(1,4,3)\}&{22}\\
\{(1,1,2),(3,4,4),(2,2,2),(4,3,4)\}&{22}\\
\{(4,3,2),(1,1,4),(2,3,3),(3,1,1)\}&{22}\\
\{(1,1,2),(4,3,2),(4,4,3),(2,1,3)\}&{24}\\
\{(1,1,1),(3,4,3),(2,2,1),(4,3,3)\}&{24}\\
\{(2,2,4),(3,4,2),(1,3,4),(4,1,2)\}&{24}\\
\{(2,4,2),(2,3,3),(1,4,1),(1,3,4)\}&{26}
\end{array}
\end{array}
\]
}
\caption{Increasing 2-ranks by repeated GM-switching in $G_+(3)$ (left) and $G'_+(3)$ (right)}
\end{table}

The upper bound for the $\rank$ of a graph with parameters $P_\pm(3)$ is 28.
So only the existence of one with $\rank$ 28 is unsolved.
If $G$ is a SRG with parameters $P_\pm(3)$ with $\rank(G)=26$, and $\1\not\in\col(G)$,
then from Lemma~\ref{Seidel2} it follows that isolating a vertex by Seidel switching gives an
SRG $G'$ with parameter set $P_0(3)$ and $\rank(G')=26$, the only open case for $P_0(3)$.
Unfortunately, it turns out that $\1\in\col(G)$ for every graph $G$ in Table~\ref{P-} and \ref{P+}.
This is not very surprising, since we know that $\1\in\col(G_\pm(3))$ and $\1\in\col(G'_\pm(3))$, and in Section~\ref{GM}
we observed that $\1$ remains in the column space of the adjacency matrix if GM-switching increases the $\rank$.

\section{SRGs with parameters $P_0(m)$ and $P_\pm(m)$}

The computer result from the previous section and the graph product introduced in Section~\ref{prod}
lead to the following result.

\begin{thm}\label{main}
{~}\\[-15pt]
\begin{description}
\item[(i)]  There exist SRGs with parameter set $P_0(m)$ and $2$-rank $r$ for every even $r\in[2m,2(m+9\lfloor\frac{m}{3}\rfloor)]$.
\item[(ii)] There exist SRGs with parameter set $P_+(m)$ and $2$-rank $r$ for every even $r\in[2(m+1),2(m+1+9\lfloor\frac{m}{3}\rfloor)]$.
\item[(iii)]There exist SRGs with parameter set $P_-(m)$ and $2$-rank $r$ for every even $r\in[2(m+1),2(m+1+9\lfloor\frac{m}{3}\rfloor)]$.
\end{description}
\end{thm}
\noindent
{\bf Proof.}
Put $\ell=\lfloor\frac{m}{3}\rfloor$, and let $G_1,\ldots,G_\ell$ be graphs coming from normalized graphical Hadamard matrices of order $64$, which are given in Table~\ref{P0}
(so $G_i$ is a SRG with parameters $P_0(3)$ extended with an isolated vertex),
and let $G_0$ be the graph of the normalized Hadamard matrix of order $4^{m-3\ell}$.
Put $r_i=\rank(G_i)$ for $i=0,\ldots,\ell$, and define $G=G_0\otimes G_1 \otimes\cdots\otimes G_\ell$.
Then $G$ is a SRG with parameters $P_0(m)$, extended with an isolated vertex,
and Theorem~\ref{prod2rank}(iii) implies that $\rank(G)=r_0+r_1+\cdots+r_\ell$.
Now by the results in the previous section, for $i=1,\ldots,\ell$, we can choose for $r_i$ any even number in $[6,24]$.
This proves $(i)$.

The proofs of (ii) and (iii) go similarly.
Let $G_1,\ldots,G_\ell$ be SRGs with parameters $P_\pm(3)$ given in Table~\ref{P-} and \ref{P+}
(so $G_i$ comes from a regular graphical Hadamard matrices of order $64$).
For $G_0$ we take $K_1$ if $m-3\ell=0$, $2K_2$ if $m-3\ell=1$, and $2K_2\H 2K_2$ if $m-3\ell=2$.
Again $G=G_0\otimes G_1 \otimes\cdots\otimes G_\ell$, and $r_i=\rank(G_i)$ for $i=0,\ldots,\ell$.
Then $G$ is a SRG with parameters $P_\pm(m)$, and for each of $r_1,\ldots,r_\ell$ we can take any even value in $[8,\ldots,26]$.
We have seen that $\1\in\col(G_i)$ for each $G_i$, unless $i=0$ and $m=3\ell$.
Therefore Theorem~\ref{prod2rank} gives $\rank(G)=r_0+r_1+\cdots+r_\ell - 2\ell$ if $m>3\ell$,
and $\rank(G)=r_1+\cdots+r_\ell - 2\ell+2$ if $m=3\ell$.
So $\rank(G)$ can become any even number in $[2(m+1),2(m+1+9\lfloor\frac{m}{3}\rfloor)]$.
If an odd number of graphs $G_1,\ldots,G_\ell$ have parameters $P_+(3)$, then $G$ has parameters $P_+(m)$,
otherwise $G$ has parameters $P_-(m)$.
\qed

By Lemma~\ref{Seidel2},
isolating by Seidel switching a vertex of a SRG $G$ with parameters $P_\pm(m)$ and 2-rank $r$,
gives a SRG with parameters $P_0(m)$ and 2-rank $r-2$ if $\1\in\col(G)$, and $r$ otherwise.
Since each graph $G$ from Tables~\ref{P-} and \ref{P+} has $\1\in\col(G)$, case (i) of Theorem~\ref{main}
can also be obtained from case (ii), or (iii).
\\

Two Hadamard matrices are {\em equivalent} if one can be obtained from the other by row and column permutation and multiplication
of rows and columns by $-1$.
Clearly each graphical Hadamard matrix is equivalent to a normalized graphical Hadamard matrix,
and by Lemma~\ref{Seidel2}, the SRGs from equivalent normalized graphical Hadamard matrices have the same 2-rank.
So case (i) of Theorem~\ref{main} gives:
\begin{cor}\label{H}
The number of nonequivalent graphical Hadamard matrices of order $4^m$ is unbounded.
\end{cor}
Lemma~\ref{Seidel2} implies that the 2-ranks of graphs from equivalent graphical Hadamard matrices differ by at most $2$.
Therefore Theorem~\ref{main} also implies that the statement of Corollary~\ref{H} remains true if we restrict to
regular graphical Hadamard matrices.

Corollary~\ref{H} may be an open door.
For several values of $m$ there exist a large number of nonequivalent (regular graphical) Hadamard matries of order $4^m$,
and by taking Kronecker products this leads to exponentially many different constructions.
However, we are not aware of another result that proves the nonequivalence of an unbounded number of these constructions.

SRGs with parameters $P_+(m)$ are known as max energy graphs, see~\cite{H}.
So Theorem~\ref{main}(ii) implies that the number of nonisomorphic max energy graphs of order $4^m$ is unbounded.

\subsection*{Acknowledgments}
Aida Abiad is supported by the \emph{Elinor Ostrom Research Grant}.
Steve Butler is supported by a grant from the Simons Foundation (\#427264).

\end{document}